\documentclass[12pt]{article} 
\usepackage{amsmath,amssymb,eucal}
\usepackage{stackrel}
\font\tenrm=cmr10

\font\cmssl=cmss10 at 12 pt  
\font\bigss=cmssdc10 scaled 2300

\font\cmsslll=cmss10 at 14 pt

\newcommand{\e}{\epsilon}

\newcommand{\s}{\sigma}

\newcommand{\bC}{\mathbb{C}}  
\newcommand{\bR}{\mathbb{R}}  
\newcommand{\bZ}{\mathbb{Z}}

\newcommand\GL{\mathrm{GL}}

\newcommand{\p}{\partial}  
    
\renewcommand{\square}{\kern1pt\vbox  
               {\hrule height 0.6pt\hbox{\vrule width 0.6pt\hskip 3pt  
    \vbox{\vskip 6pt}\hskip 3pt\vrule width 0.6pt}\hrule height0.6pt}  
    \kern1pt}  

\newcommand{\ra}{\rightarrow}

\DeclareMathOperator\tr{tr\;}

\newcommand{\wt}{\widetilde}

\newtheorem{Pb}{Problem}

\newtheorem{Th}{Theorem}  
\newtheorem{Prop}{Proposition}  
  
\newtheorem{Cor}{Corollary}  
\newtheorem{Lem}{Lemma}  
\newtheorem{Def}{Definition}
\newtheorem{Rem}{Remark}
\newcommand{\bP}{\begin{Pb}\ \ } 
\newcommand{\eP}{\end{Pb}}  
\newcommand{\bt}{\begin{Th}\ \ }  
\newcommand{\et}{\end{Th}}  
\newcommand{\bp}{\begin{Prop}\ \ }  
\newcommand{\ep}{\end{Prop}}  
\newcommand{\bc}{\begin{Cor}\ \ }  
\newcommand{\ec}{\end{Cor}}  
\newcommand{\bl}{\begin{Lem}\ \ }  
\newcommand{\el}{\end{Lem}}  
\newcommand{\bd}{\begin{Def}\ \ }  
\newcommand{\ed}{\end{Def}}  
\newcommand{\br}{\begin{Rem}\ \ }  
\newcommand{\er}{\end{Rem}}  
  
\newcommand{\pf}{\noindent{\it Proof:\ \ }}  
\newcommand{\qed}{\hfill\square}

\newcommand{\be}{\begin{equation}}  
\newcommand{\ee}{\end{equation}}  
  
\newcommand\re[1]{(\ref{#1})}  
\newcommand{\arr}{\begin{array}{rlll}}  
\newcommand{\ea}{\end{array}}  
\newcommand{\bea}{\begin{eqnarray}}  
\newcommand{\eea}{\end{eqnarray}}  
\newcommand{\bean}{\begin{eqnarray*}}  
\newcommand{\eean}{\end{eqnarray*}}  
  
\catcode`@=11  
\@addtoreset{equation}{section}  
\catcode`@=12

\begin{document}  
%\begin{titlepage}
 \rightline{} 
%\rightline{hep-th/yymmnnn}  
%\rightline{draft: \today}  
\vskip 1.5 true cm  
\begin{center}  
{\bigss  Classification of complete projective special real surfaces}\\[.5em]
\vskip 1.0 true cm   
{\cmsslll  V.\ Cort\'es$^1$, M.\ Dyckmanns$^2$ and D.\ Lindemann$^2$} \\[3pt] 
$^1${\tenrm   Department of Mathematics\\  
and Center for Mathematical Physics\\ 
University of Hamburg\\ 
Bundesstra{\ss}e 55, 
D-20146 Hamburg, Germany\\  
cortes@math.uni-hamburg.de}\\[1em]  
$^2${\tenrm Department of Mathematics\\  
University of Hamburg\\ 
Bundesstra{\ss}e 55, 
D-20146 Hamburg, Germany}\\ [1em]
February 18, 2013 
\end{center}  
\vskip 1.0 true cm  
%%%%%%%%%%%%%%%%%%%%%%%%%%%%%%%%%%%%%%%%%%%%%%%%%%%%%%%%  
\baselineskip=18pt  
\begin{abstract}  
\noindent  
We determine all complete projective special
real surfaces. By the supergravity r-map, they give rise to complete 
projective special K\"ahler manifolds of dimension 6,
which are distinguished by the image of their 
scalar curvature function.  By the supergravity c-map, the latter manifolds 
define in turn complete quaternionic 
K\"ahler manifolds of dimension 16.
\end{abstract}

%\end{titlepage}  
\tableofcontents
\section*{Introduction}
Projective special real manifolds 
first occurred as the scalar manifolds of certain supergravity theories in 
five space-time dimensions \cite{GST,DV1}, see Definition \ref{psrDef} below. 
Their geometry is encoded in a 
homogeneous cubic polynomial. A typical example occurring in string theory is 
the geometry defined by the cubic form on $H^{1,1}(X,\mathbb{R})$ for a 
K\"ahler manifold $X$ of complex dimension 3. 
It was shown in \cite{CHM}, using constructions from supergravity, that any 
complete projective special real manifold of  dimension $n$ defines a complete
projective special K\"ahler manifold of (real) dimension  $2n+2$ and 
a complete quaternionic K\"ahler manifold of dimension $4n+8$. 
For that reason it is interesting 
to find examples of complete projective special real manifolds.   
Let us also mention that the completeness (or incompleteness)  
of the scalar manifold in the underlying supergravity theories is 
related to the global behaviour of solutions to the equations of motion. 
This is due to the fact that the scalar fields 
of the theory cannot approach infinity along a trajectory of finite 
length if the manifold is complete. 

In this paper we classify all complete projective special real surfaces:

\bt \label{mainThm} There exist precisely five discrete examples and a one-parameter family of complete projective special real surfaces, 
up to isomorphism: 
\begin{itemize}
\item[a)] $\{ (x,y,z)\in \bR^3| xyz=1,\; x>0,\; y>0\}$, 
\item[b)] $\{ (x,y,z)\in \bR^3| x(xy-z^2)=1,\; x>0\}$, 
\item[c)] $\{(x,y,z)\in \bR^3| x(yz+x^2)=1,\; x<0,\; y>0\}$, 
\item[d)] $\{(x,y,z)\in \bR^3| z(x^2+y^2-z^2)=1,\; z<0\}$,   
\item[e)] $\{ (x,y,z)\in \bR^3| x(y^2-z^2) + y^3=1,\;  y<0,~x>0\}$,
\item[f)] $\{ (x,y,z)\in \bR^3| y^2z-4x^3+3xz^2+bz^3=1,\; z<0,\; 2x>z\}$, where $b\in(-1,1)\subset \mathbb{R}$.
\end{itemize}
\et 

In Sections \ref{StatementSec}--\ref{ClassSec} we prove Theorem \ref{mainThm}  
essentially by first determining all homogeneous cubic polynomials $h$ on $\bR^3$ such that
the surface $\{ h=1\}\subset \bR^3$ has a strictly locally convex component
$\mathcal{H}$ of hyperbolic type. It turns out that, up to linear transformations,
the resulting surfaces $\mathcal{H}$  of hyperbolic type are precisely those listed in Theorem \ref{mainThm}. 
Then we prove in all cases that $\mathcal{H}$ 
is complete with respect to the Riemannian metric  $g_\mathcal{H}$ induced by the 
Hessian of $-h$.  At present we do not know, in dimensions $n\ge 3$, whether a projective special
real manifold $\mathcal{H}\subset \bR^{n+1}$ which is closed as a subset of $\bR^{n+1}$
is necessarily complete with respect to the metric $g_\mathcal{H}$, see open problem on page
\pageref{openPb}. The converse statement, however, can be easily proven in all dimensions. 

In Section \ref{Curv_genSec} we provide general formulas for the  curvature  of 
the K\"ahler manifolds obtained from the generalized r-map  defined
in \cite{CHM}. This is applied in Section \ref{lastSec} to the projective 
special K\"ahler manifolds of complex dimension 3 obtained from the projective special 
real surfaces of Theorem \ref{mainThm} . Computing the scalar curvature we prove, in particular, that
the complete projective special real manifolds c)--f) of Theorem \ref{mainThm} as well as the corresponding complete projective 
special K\"ahler manifolds are not locally homogeneous as Riemannian manifolds. 
The family of projective special K\"ahler manifolds associated with the 
Weierstra{\ss} polynomials  $h_b= y^2z-4x^3+3xz^2+bz^3$, $b\in (-1,1)$, seems to be the first example of a continuous family
of complete projective special K\"ahler manifolds. 

Let us finally mention that applying the supergravity 
c-map to our examples one obtains complete quaternionic K\"ahler manifolds of dimension
16, which have negative scalar curvature and cohomogeneity less than or equal to $2$. 
More precisely, there exists a group of isometries 
which has cohomogeneity $k$, where $k$ is the cohomogeneity of the initial projective special real
manifold under the full group of linear automorphisms, that is $k=0$ in the cases
a) and b), $k=1$ in the case c) and d) and $k=2$ in the cases e) and f) . 
We expect that the same is true for the full isometry group of the quaternionic K\"ahler manifold. 
%The detailed study of these quaternionic K\"ahler manifolds is left for 
%future investigation. 

 \noindent
{\bf Acknowledgments} 
This work is part of a research 
project within the RTG 1670 ``Mathematics inspired by String Theory'',
funded by the Deutsche Forschungsgemeinschaft (DFG).

\section{Statement of the classification problem}\label{sectionStatement}
\label{StatementSec}
Let $h$ be a homogeneous cubic
polynomial function on $\bR^{n+1}$ and $U\subset \bR^{n+1}$ 
a domain invariant under multiplication 
by positive numbers such that $h|_U>0$. 
Then 
\[ \mathcal{H} :=\{ x\in U| h(x)=1\}\subset U\]
is a smooth hypersurface and 
$-\partial^2h$ induces a symmetric tensor
field $g_{\mathcal{H}}$ on $\mathcal{H}$.

\bd \label{psrDef} If $g_{\mathcal{H}}$ is positive  
definite then the hypersurface $\mathcal{H}\subset \bR^{n+1}$ 
is called a {\cmssl projective special real manifold}. 
Two projective special real manifolds $\mathcal{H}, \mathcal{H}'
\subset \bR^{n+1}$ are called {\cmssl isomorphic} if there exists 
$\varphi \in \GL(n+1)$ 
mapping  $\mathcal{H}$ to $\mathcal{H}'$.  
\ed 

\br
\rm{\\1) In the above situation $\varphi|_{\mathcal{H}} : 
\mathcal{H}\ra \mathcal{H}'$ is an isometry.\\
2) For any projective special real manifold $\mathcal{H}
\subset \bR^{n+1}$ the tensor field 
$-\partial^2h$ is a Lorentzian metric on $U=\bR^{>0}\cdot 
\mathcal{H}\subset \bR^{n+1}$. 
In particular, $(-1)^n\det \partial^2h>0$ on $U$.}
\er

The classification of complete projective special real surfaces
up to isomorphism can be separated into two problems. 
First we have to determine all homogeneous cubic polynomials
$h$ on $\bR^3$, up to linear transformations, which are {\cmssl hyperbolic}, 
that is admit a point $p\in \bR^3$ 
such that $h(p)>0$ and $\partial^2h$ is negative definite on the 
kernel of $dh_p$.  We can assume that $h(p)=1$. 
Then there exists a maximal connected 
neighborhood $\mathcal{H}$ of $p$ in the level set $\{ h=1\}$ such that  
$\mathcal{H}$ is a  projective special real surface. The second problem is
then to check whether $(\mathcal{H},g_{\mathcal{H}})$ is a complete Riemannian 
manifold.

\section{Classification of cubic polynomials}\label{sectionClassPolynomials}
In this section we provide the needed classification of homogeneous real cubic 
polynomials $h$ in three variables up to linear transformations. We will
say that two polynomials are {\cmssl equivalent} if they are related 
by a linear transformation. This problem 
is equivalent to the classification of cubic curves 
in the real projective plane.  The study of real plane cubic curves 
goes back to Newton \cite{N}. For the classification of 
\emph{complex} plane cubic curves up to projective transformations 
see the textbooks \cite{Brieskorn,Hulek}.    

Let us first consider the case when $h$ is reducible, that is
a product of homogeneous polynomials of degree 1 or 2. 

\bp\label{redProp} Any reducible homogeneous cubic polynomial on $\bR^3$ is 
equivalent to one of the following:
\begin{enumerate}
\item[(i)] $x^3$,
\item[(ii)] $x^2y$,
\item[(iii)] $xy(x+y)$,
\item[(iv)] $xyz$,
\item[(v)] $x(x^2+y^2)$,
\item[(vi)] $z(x^2+y^2)$,
\item[(vii)] $x(x^2+y^2+z^2)$,
\item[(viii)]  $z(x^2+y^2-z^2)$,
\item[(ix)]  $x(x^2+y^2-z^2)$,
\item[(x)]  $(y+z)(x^2+y^2-z^2)$.
\end{enumerate}
\ep 

The proof is a simple exercise. Notice that the first 4 cases are products
of linear factors and are the same as in the complex case. The remaining
6 polynomials contain an irreducible quadratic polynomial 
as a factor.  Over $\bC$ there are only 2 such polynomials.   

Next we consider the case of singular curves $C=\{ h=0\} \subset \bR P^2$. 
\bp\label{singProp} Any irreducible homogeneous cubic polynomial $h$ on $\bR^3$ such that
the curve $C=\{ h=0\} \subset \bR P^2$ has a singularity is equivalent 
to one of the following:
\begin{enumerate}
\item[(xi)] $x(y^2+z^2)+y^3$,   
\item[(xii)] $x(y^2-z^2)+y^3$,
\item[(xiii)] $xz^2+y^3$. 
\end{enumerate}
\ep 

\pf  We can assume that $h$ and $dh$  
vanish at $P=(1,0,0)$. Then we decompose $h=xq+r$, 
where $q=q(y,z)$, and $r=r(x,y,z)\neq 0$ does not contain
any monomial summands linear in $x$. Now the conditions
$h(P)=0$ and $\partial_yh(P)=\partial_zh(P)=0$ easily imply that 
$r=r(y,z)$. By a linear transformation we can obviously assume that 
$q\in \{ y^2+z^2, y^2-z^2,z^2\}$.  In the last two cases, one can
use the same linear transformation as in the complex case
to bring $h$ to the form (xii) and (xiii), respectively, 
cf.\ \cite{Hulek}.  Therefore it suffices to consider the case
$q=y^2+z^2$. The vector 
space $S\cong S^3(\bR^2)^*$ of homogeneous cubic polynomials  
in the variables $(y,z)$ is decomposed as a sum of two irreducible
$\mathrm{O}(2)$-modules: 
\[ S^3(\bR^2)^* = (\bR^2)^*q \oplus \mathrm{span}\{ y^3-3yz^2,z^3-3y^2z\}.\] 
Using a homothety in the $(y,z)$-plane we can thus assume
that $r\equiv y^3-3yz^2$ modulo $(\bR^2)^*q$. 
Furthermore, by a linear transformation preserving $y$ and $z$ we can 
freely change the $(\bR^2)^*q$-component of $r$. For instance,
we can take it to be $3yq$, which implies $r=4y^3$. Now 
it suffices to rescale $x$ and $(y,z)$ to bring $h$ to the 
form (xi). 
\qed 

Finally, the classification 
of smooth irreducible real cubic curves is provided 
by the Weierstra{\ss} normal form. 

\bp\label{propSmoothWCubics} Let $h$ be an 
irreducible homogeneous cubic polynomial on $\bR^3$ such that
the curve $C=\{ h=0\} \subset \bR P^2$ is smooth. Then $h$ is 
equivalent to a Weierstra{\ss} cubic polynomial 
\[ y^2z-4x^3+axz^2+bz^3,\]
of nonzero discriminant $a^3-27b^2$, for some $a, b\in \bR$. 
\ep  

\pf 
By Bezout's theorem, we know that $C$ has 9 complex 
inflection points \cite{Brieskorn,Hulek}. Since, the imaginary 
inflection points occur in pairs one of them has to be real. 
Therefore the same proof as in the complex case applies \cite{Hulek}. 
\qed 

\br \rm{Note that some of the Weierstra\ss~cubic polynomials given in Proposition \ref{propSmoothWCubics} are linearly equivalent. From the classification of smooth cubics over $\mathbb{C}$, we know that two Weierstra{\ss} cubics are inequivalent if they have different j-invariants (see \cite{Brieskorn,Hulek})}
\[j(a,b):=\frac{a^3}{a^3-27b^2}.\]
\er

\section{Classification of hyperbolic polynomials}
In this section we study the hyperbolicity of the polynomials given in Propositions \ref{redProp} and \ref{singProp}. The study of Weierstra{\ss} polynomials with nonzero discriminant is postponed to Subsection \ref{weierSection}.

\bp \label{hypProp} 
Let $h$ be a hyperbolic homogeneous cubic polynomial on $\bR^3$.
Then either $h$ is equivalent to a Weierstra{\ss} cubic polynomial with
nonzero discriminant or to one of the following:
\begin{enumerate}
\item[1)] $xyz$,
%a  
\item[2)]  $z(x^2+y^2-z^2)$,
%d
\item[3)]   $x(x^2+y^2-z^2)$,
%c
\item[4)]   $(y+z)(x^2+y^2-z^2)$,
%b
\item[5)]  $x(y^2+z^2)+y^3$,   
\item[6)]  $x(y^2-z^2)+y^3$,
%e
\item[7)]  $xz^2+y^3$.
\end{enumerate}
\ep 

\pf In the following we will denote by $g$ the symmetric tensor field
on the surface $\{ h=1\}$, which is induced by $-\p^2h$. We have to decide
whether $-\p^2h(p)$ is Lorentzian for some $p\in \{ h=1\}$ or, equivalently, 
whether $g_p$ is positive definite for some $p\in \{ h=1\}$. 

The polynomials (i-iii) and (v) in Proposition \ref{redProp} 
do not depend on $z$. 
Hence, their Hessian is everywhere degenerate. We claim that
also (vi) and (vii) are not hyperbolic, which leaves us with 
the cases 1)-7). In the case (vi), $\det \p^2 h = -8h$ is negative
on $\{ h=1\}$. Therefore, $h$ is not hyperbolic. 
In the case (vii),  $\det \p^2 h = 8(4x^3-h)$, 
which is positive 
precisely on those points of $\{ h=1\}$ for which 
$x>\frac{1}{\sqrt[3]{4}}$.
Next we observe that $\det \begin{pmatrix} h_{yy} & h_{yz} \\ h_{zy} & h_{zz} \end{pmatrix}=4x^2$ and $h_{zz}=2x$ are positive for $x>\frac{1}{\sqrt[3]{4}}$, where the subscripts denote partial derivatives. This shows that $\p^2 h$ is positive definite
on the set $\{ \det \p^2 h>0\}\cap \{ h>0\}$.
%Next we observe that the determinant 
%of the two-by-two leading principal minor of $\p^2 h$ is 
%a positive multiple of $3x^2-y^2$. On $\{ h=1\}$, we have 
%$y^2 = \frac{1}x-x^2 -z^2 \le  \frac{1}x-x^2$. Thus 
%$3x^2-y^2\ge 4x^2 -\frac{1}{x}$, which is positive for  
%$x>\frac{1}{\sqrt[3]{4}}$. Since $\partial^2_xh=6x$ is positive 
%on $\{ h=1\}$, this shows that $\p^2 h$ is positive definite
%on the set $\{ \det \p^2 h>0\}\cap \{ h>0\}$.
In particular,
$h$ is not hyperbolic. Now we prove that the remaining 
polynomials are all hyperbolic.  The polynomials 1), 3) and 4) 
give rise to the complete projective special real 
surfaces a)-c) in Theorem \ref{mainThm}, which were 
already discussed in \cite{CHM}.   In fact, $x(x^2+y^2-z^2)$
is equivalent to $xyz + x^3$ and $(y+z)(x^2+y^2-z^2)$
to $x(xy -z^2)$.\\
In case 2),  $\det \p^2 h=-8(4z^3+h)$, 
which is positive precisely on those points of $\{ h=1\}$ for which  
$z<-\frac{1}{\sqrt[3]{4}}$. 
Since  $\partial^2_xh=2z$, we see that the subset of 
$\{ h=1\}$ on which $g$ is Riemannian 
is nonempty and coincides with the subset on 
which $\det \p^2 h>0$. 
This subset is precisely the surface d) in Theorem \ref{mainThm}.
In fact, $z<0$ and $h=1$ imply $z\le -1<-\frac{1}{\sqrt[3]{4}}$.\\ 
In case 5),   the surface $\{ h=1\}$ is a graph 
\be  x=\frac{1-y^3}{y^2+z^2} \nonumber\ee   
over the domain $\bR^2\setminus \{ 0\}$ in the $(y,z)$-plane. 
The (nonempty) subset
on which $g$ is Riemannian is 
\be \{ (y,z)\in \bR^2|y(y^2-3z^2)-1>0\}.\nonumber\ee
This follows from $\det\partial^2h= 8(y(y^2-3z^2)-h)$, since $\p^2_xh=0$.\\    
In case 6), the surface $\{ h=1\}$ is a union of the graph
\be  \{ x=\frac{1-y^3}{y^2-z^2},\; y^2-z^2\neq 0\} \nonumber\ee
and the two vertical lines $\{ y=|z|=1\}$. Since 
$\det\partial^2h= 8(h -y(y^2+3z^2))$ and $\p^2_xh=0$, 
we see that 
the (nonempty) subset of $\{ h=1\}$ on which $-\p^2 h$ induces a Riemannian 
metric on the surface is precisely 
\be  \{ y(y^2+3z^2)<1,\; x=\frac{1-y^3}{y^2-z^2},\; y^2-z^2\neq 0 \}.\nonumber\ee 
In case 7),  $\{ h=1\}$ is a graph
\be y=\sqrt[3]{1-xz^2}\nonumber\ee 
over the $(x,z)$-plane and $g$ is Riemannian precisely on the (nonempty) 
subset \be xz^2>1.\nonumber\ee
This follows from $\det \partial^2h=-24yz^2$, since $\p^2_xh=0$.\\
\qed

\section{Classification of complete surfaces} \label{ClassSec} 
In this section we study the completeness of the 
maximal connected projective special real surfaces associated 
with the hyperbolic cubic polynomials $h$ 
described in Proposition \ref{hypProp}. 
These are precisely the connected components
of the hypersurface
\[ \mathcal{H}(h) := \{ x\in \bR^3| h(x)=1\quad
\mbox{and}\quad g_x>0\},\]
where we recall that $g_x$ is the restriction of 
the symmetric bilinear form $-\partial^2h(x)$ to
the plane $\ker dh(x)$. 

In the next theorem we determine all the complete  
and incomplete components of $\mathcal{H}(h)$ 
for the polynomials  1)-7) of 
Proposition \ref{hypProp}, up to equivalence.
The five cases which admit a complete component 
are listed first. They correspond to the surfaces a)-e) in 
Theorem \ref{mainThm}. The case of Weierstra{\ss} polynomials with nonzero discriminant will 
be analysed in the next subsection. It will lead to the family of surfaces f) in Theorem \ref{mainThm}.
\bt\label{theoremComplComp} \begin{enumerate}
\item[1)] For $h=xyz$, $\mathcal{H}(h)=\{ h=1\}$ 
has four isomorphic components, 
each of which is complete. 
%a)
\item[2)] For $h= x(xy-z^2)$,  $\mathcal{H}(h)= 
\{ (x,y,z)\in \bR^3| h(x,y,z)=1,\; x>0\}$ is connected
and complete. 
%b)
\item[3)]  For $h= x(yz+x^2)$,  $\mathcal{H}(h)=\{ (x,y,z)\in \bR^3| 
h(x,y,z)=1,\; x<\frac{1}{\sqrt[3]{4}}\}$ has  four components; 
a pair of isomorphic complete components and a pair of 
isomorphic incomplete components.  
%c)
\item[4)]  For $h= z(x^2+y^2-z^2)$, $\mathcal{H}(h)=\{ (x,y,z)\in \bR^3| 
h(x,y,z)=1,\; z<0\}$ is connected and complete. 
%d)
\item[5)]  For $h= x(y^2-z^2) + y^3$, 
\[ \mathcal{H}(h)=\{ (x,y,z)\in \bR^3| 
h(x,y,z)=1,\; y^2-z^2\neq 0,\; y(y^2+3z^2)<1\}\] 
has 
four components; a complete component, 
a pair of isomorphic incomplete components and a further incomplete
component.  
%e)
\item[6)] For $h=x(y^2+z^2)+y^3$, $\mathcal{H}(h)=\{(x,y,z)\in\mathbb{R}^3|h(x,y,z)=1,~y(y^2-3z^2)>1\}$ has three components; a pair of isomorphic incomplete components and a further incomplete component.
\item[7)] For $h=xz^2+y^3$, $\mathcal{H}(h)=\{(x,y,z)\in\mathbb{R}^3|h(x,y,z)=1,~y<0\}$ has two components. They are isomorphic and incomplete.
\end{enumerate}
\et 

\pf 
1) The group $\bZ^2_2$ acts 
by $(x,y,z) \mapsto (\e_1 x,\e_2y, \e_1\e_2z)$,  $\e_1 ,\e_2 \in \{ \pm 1\}$,  
on the level set $\{ h=1\}$  
permuting its four components. Therefore, the statement follows
from the fact \cite{CHM} that the tensor field $g$ is positive definite 
and complete on the component $\{ (x,y,z)\in \bR^3| xyz=1,\; x>0,\; y>0\}$.\\
2) The description of  $\mathcal{H}(h)$ follows from 
$\det \partial^2h = 8x^3$ and $\partial^2_yh=0$, the 
completeness from \cite{CHM}.\\  
3) The description of  $\mathcal{H}(h)$ follows from 
$\det \partial^2h = 2(h-4x^3)$ and $\partial^2_yh=0$.
Notice that $\mathcal{H}(h)$ is a graph over the union 
of the following four domains in the $(x,y)$-plane:
$\{ x<0, y>0\},\{ x<0, y<0\}, \{ 0<x<\frac{1}{\sqrt[3]{4}}, y>0\} ,  
\{ 0<x<\frac{1}{\sqrt[3]{4}}, y<0\}$. The corresponding 
components of $\mathcal{H}(h)$ are related by the 
involution $(x,y,z) \mapsto (x,-y,-z)$. So up to isomorphism, 
it suffices to consider the two components 
\begin{eqnarray} \{ (x,y,z)\in \bR^3| 
h(x,y,z)=1,\; x< 0, y>0\}, \nonumber\\ 
\{ (x,y,z)\in \bR^3| 
h(x,y,z)=1,\;  0<x<\frac{1}{\sqrt[3]{4}}, y>0\}.\nonumber
\end{eqnarray}       
The first one is complete by \cite{CHM} and the second one 
is incomplete. This follows from the fact that the second component
has nonempty boundary. The boundary is given by the curve      
\[ \left\{ \left.
\left( \frac{1}{4^{1/3}},y,\frac{3}{4^{2/3}y}\right)\right| 
y> 0 \right\} . \]   
4) The description of  $\mathcal{H} = \mathcal{H}(h)$ 
was obtained in the proof of 
Proposition \ref{hypProp}. In order to prove the completeness
let us first remark that $\mathcal{H}(h)$ is a surface of revolution. 
More precisely, it is the graph of the function 
\[ (x,y)\mapsto z=\varphi (\rho ),\]
where $\rho = r^2 = x^2 +y^2$ and 
\[ \varphi : [0,\infty ) \ra 
(-\infty , -1]\] is the inverse of the strictly decreasing
function 
\[ f : (-\infty , -1]\ra [0,\infty ),\quad 
z\mapsto \rho = f(z)=\frac{1}{z}+z^2.\]   
Let us first calculate the metric $g=-\partial^2h|_\mathcal{H}$ in the 
coordinates $(x,y)$. Using that $z=\varphi (\rho )$, we obtain 
\[ \frac{1}{2}g= -z(dx^2 + dy^2) +3zdz^2- 2(xdx+ydy)dz
=-\varphi (\rho )(dx^2+dy^2) +\varphi'(\rho )(3\varphi (\rho )\varphi'(\rho )- 
1)d\rho^2  . \]
Rewriting $dx^2+dy^2= dr^2+r^2ds^2= \frac{1}{4\rho}d\rho^2+\rho ds^2$  
in polar coordinates $(r,s)$ in the $(x,y)$-plane  
we arrive at 
\be \label{f1f2Equ} g|_{\{ \rho >0\} }= 
2f_1(\rho )d\rho^2 +2f_2(\rho )ds^2 , \ee
where 
\begin{eqnarray*}
f_1(\rho ) &=& -\frac{1}{4\rho}\varphi (\rho )+\varphi'(\rho )(3\varphi 
(\rho )\varphi'(\rho )- 
1)\\
f_2(\rho ) &=& -\varphi (\rho )\rho . 
\end{eqnarray*}
The metric \re{f1f2Equ} is of the type considered in 
Section 1 of \cite{CHM} (cf.\ Lemma 
 \ref{lemmaCompleteness} below). Therefore, for all $a>0$, the completeness of
$\mathcal{H}$ is equivalent to 
\be \label{intEqu} \int_a^\infty \sqrt{f_1(\rho )}d\rho = \infty.\ee 
A straightforward calculation shows the following asymptotics when 
$\rho \ra \infty$: 
\[ f_1(\rho)  = \frac{3}{4\rho^2} + O(\rho^{-7/2}),\]
which implies \re{intEqu}.\\
5) The description of  $\mathcal{H}(h)$ was obtained in the proof of 
Proposition \ref{hypProp}. The surface is a graph over
the union of the following four domains in the $(y,z)$-plane: 
\[ \{ y<0, |z|<|y|\},\; \{ 0<y<1, |z|< \mathrm{min}(y, f(y))\},\; \{ \e z> 0, |y| < 
|z|, y< 
f^{-1}(|z|)\},\]
where $\e=\pm 1$, $f : (0,1) \stackrel{\sim}{\ra} (0,\infty )$ is the strictly 
decreasing function
$f(y)= \frac{1}{\sqrt{3}}(\frac{1}{y}-y^2)^{1/2}$ and $f^{-1}: (0,\infty )
 \stackrel{\sim}{\ra} (0,1)$ denotes its inverse. The involution 
$(x,y,z) \mapsto (x,y,-z)$ acts on $\mathcal{H}(h)$ preserving
the first two components and interchanging the last two. The last three
components are incomplete as they have a nonempty boundary. 
The boundary of the second component is 
\[ \{ (x,y,z)\in \bR^3| y_0< y\le1,\; |z|=f(y),\; 
x= \frac{1-y^3}{y^2-z^2}\},\]
where $y_0$ is the unique fixed point of $f$. 
The boundary of the third one ($\epsilon=1$) is 
\[ \{ (x,y,z)\in \bR^3| 0< y<y_0,\; z=f(y),\; 
x= \frac{1-y^3}{y^2-z^2}\}.\] 
Now we show that the first component, namely  
\[ \mathcal{H} 
:= \{ (x,y,z)\in \bR^3|y<0, |z|<|y|, x= 
\frac{1-y^3}{y^2-z^2}\},\]
is complete. 
The metric $g$ is given by 
\[ \frac{1}{2}g = -(x+3y)dy^2+xdz^2-2dx(ydy-zdz),\]
where $(y,z)$ is restricted to the domain
$\{ y<0, |z|<|y|\}$ and $x= \frac{1-y^3}{y^2-z^2}$.
Using the function $s= y^2-z^2>0$, we rewrite 
this as 
\[ \frac{1}{2}g = -3ydy^2-x(dy^2-dz^2)-dxds.\]
Eliminating $x=\frac{1-y^3}{s}$ and 
\[dx=-\frac{3y^2dy}{s}-\frac{(1-y^3)ds}{s^2},\] 
we get 
\[ \frac{1}{2}g = -3ydy^2-(1-y^3)\frac{dy^2-dz^2}{s}+3y^2dyd\s  + 
(1-y^3)d\s^2,\]
where $\s =\ln s$. 
Using the coordinates 
\[ t_\pm := \ln (|y|\pm z),\]
such that 
\[ \s  = t_++t_-,\]
this is 
\[ \frac{1}{2}g = -3ydy^2-(1-y^3)dt_+dt_-+3y^2dy(dt_++dt_-)  + 
(1-y^3)(dt_++dt_-)^2.\]
Notice that the coordinates $(t_+,t_-)$ define a diffeomorphism
$\mathcal{H}\cong \bR^2$. 
Eliminating 
\[ y= -\frac{e^{t_+}+e^{t_-}}{2},\]
we get 
\begin{eqnarray*} \frac{1}{2}g &=& 
(1+\frac{1}{8}(e^{3t_+} + e^{3t_-}))(dt_+^2+dt_-^2)
+ (1-\frac{1}{4}(e^{3t_+} + e^{3t_-}))dt_+dt_-\\
&=&  dt_+^2+dt_-^2 + dt_+dt_- + \frac{1}{8}(e^{3t_+} + e^{3t_-}))(dt_+-dt_-)^2\\
&\ge& dt_+^2+dt_-^2 + dt_+dt_- \ge \frac{1}{2}(dt_+^2+dt_-^2). 
\end{eqnarray*}
So $g$ is bounded from below by the complete metric $dt_+^2+dt_-^2$ and, hence,
is itself complete.\\
6) The description of  $\mathcal{H}(h)$ follows from $\det\partial^2 h=8(y(y^2-3z^2)-h)$ and $\partial^2_xh=0$. $\mathcal{H}(h)$ is a graph $x=\frac{1-y^3}{y^2+z^2}$ over the union of the following three domains in the $(y,z)$-plane:
$\{y<0,~3z^2>y^2-\frac{1}{y},~z>0\}$, $\{y<0,~3z^2>y^2-\frac{1}{y},~z<0\}$, $\{y>0,~3z^2<y^2-\frac{1}{y}\}$. The first two correspond to components of $\mathcal{H}(h)$ that are related by the involution $(x,y,z)\mapsto (x,y,-z)$. $(y,z)=(-1,\sqrt{2/3})$ and $(y,z)=(1,0)$ are points in the boundary of the first and third domain respectively. Since $|x|=|\frac{1-y^3}{y^2+z^2}|<\infty$ for $(y,z)\neq (0,0)$, the corresponding components of $\mathcal{H}(h)$ have nonempty boundary. Hence, all three components of $\mathcal{H}(h)$ are incomplete. \\
7) The description of $\mathcal{H}(h)$ follows from $\det\partial^2 h=-24yz^2$ and $\partial^2_xh=0$. $\mathcal{H}(h)$ is a graph $y=\sqrt[3]{1-xz^2}$ over the union of the following two domains in the $(x,z)$-plane: $\{xz^2>1,z>0\}$, $\{xz^2>1,~z<0\}$. The corresponding components of $\mathcal{H}(h)$ are related by the involution $(x,y,z)\mapsto (x,y,-z)$. $(x,y,z)=(1,0,1)$ is a point in the boundary of the first component. Hence both components of $\mathcal{H}(h)$ are incomplete.
\qed
  
\pf (of Theorem \ref{mainThm})
In Section \ref{sectionClassPolynomials}, we classified all homogeneous cubic polynomials, up to linear equivalence. They fall into three classes: reducible polynomials, irreducible polynomials $h$ such that $\{h=0\}\subset\mathbb{R}P^2$ is a singular curve and irreducible polynomials such that the corresponding cubic curve in $\mathbb{R}P^2$ is smooth. For the first two classes, all hyperbolic polynomials, i.e. all polynomials that define a projective special real surface are classified in Proposition \ref{hypProp}. Theorem \ref{theoremComplComp} then classifies all complete projective special real surfaces defined by the polynomials given in Proposition \ref{hypProp}. This gives the surfaces a)-e). The polynomials in the third class are all equivalent to Weierstra{\ss} polynomials of nonzero discriminant (Proposition \ref{propSmoothWCubics}). They are studied in the next subsection. According to Corollary \ref{corWei}, the surfaces in the one-parameter family f) are, up to equivalence, the only closed projective special real surfaces defined by Weierstra{\ss} polynomials of nonzero discriminant. Proposition \ref{propWeiCompl} in combination with Lemma \ref{lemmaWeiEquiv} shows that all surfaces in the family f) are complete.
\qed

\br \rm{More precisely, we have classified all projective special real surfaces that are closed in $\mathbb{R}^3$ and we have shown that all closed projective special real surfaces are complete.}
\er

%Our classification of complete projective special real
%surfaces, see Theorem \ref{mainThm}, together with the classification 
%of complete projective special real curves, 
%see Corollary 4 of \cite{CHM},  
%imply the following corollary.  

Together with the classification of all closed and of all complete projective special real curves in \cite{CHM}, we obtain the following corollary:

\bc\label{corClosedCom}A projective special real manifold $\mathcal{H}\subset \bR^{n+1}$ 
of dimension $n\le 2$ is complete if and only if  $\mathcal{H}$ is a 
closed subset of $\bR^{n+1}$.  
\ec 
\underline{\textbf{Open problem:}}\emph{
Does the statement of Corollary \ref{corClosedCom} hold in all dimensions?\label{openPb}}

\subsection{Complete surfaces defined by Weierstra{\ss} polynomials with nonzero discriminant}\label{weierSection}
In this subsection, we will study Weierstra{\ss} polynomials
\[h^{(a,b)}:=y^2z-4x^3+axz^2+bz^3\]
with $a^3-27b^2\neq 0$ and show that for positive discriminant, they define a one-parameter family of complete projective special real surfaces and that for negative discriminant, all connected components of $\{x\in\mathbb{R}| h(x)=1,~g_x>0\}$ are incomplete.

First, we study the connected components of $\{h=1\}$. In the case of positive discriminant, we can restrict ourselves to Weierstra\ss~cubics with $a=3$, according to the following lemma.

\bl\label{lemmaEquivalencePolynomial}
Let $h^{(a,b)}=y^2z-4x^3+axz^2+bz^3$ be a Weierstra{\ss} cubic polynomial with positive discriminant $a^3-27b^2>0$. Then $h^{(a,b)}$ is linearly equivalent to $h^{(3,\wt{b})}$ with $-1<\wt{b}<1$.
\el
\pf
$a^3-27b^2>0$ implies $a>0$. Defining $\wt{x}:=x$, $\wt{y}:=\left(\frac{3}{a}\right)^{\frac{1}{4}}y$, $\wt{z}:=\sqrt{\frac{a}{3}}\, z$ and $\wt{b}:=\left(\frac{3}{a}\right)^{\frac{3}{2}}b$, we obtain  $h^{(a,b)}(x,y,z)=\tilde{y}^2\tilde{z}-4\tilde{x}^3+3\tilde{x}\tilde{z}^2+\tilde{b}\tilde{z}^3=h^{(3,\wt{b})}(\wt{x},\wt{y},\wt{z})$ and
$a^3-27b^2=a^3(1-\wt{b}^2)>0$ if and only if $-1<\wt{b}<1.$
\qed

\bp\label{propConnComponents}
Let $h=y^2z-4x^3+axz^2+bz^3$ with $a,b\in \mathbb{R}$ such that the discriminant $a^3-27b^2$ is nonzero. Then $\{(x,y,z)\in\mathbb{R}^3|h(x,y,z)=1\}$ has
\begin{enumerate}
\item[a)] one connected component for $a^3-27b^2<0$ and
\item[b)] two connected components for $a^3-27b^2>0$. For $a=3$ (and hence $-1<b<1$), one of them is given by $\mathcal{H}^{(3,b)}:=\{(x,y,z)\in\mathbb{R}^3| h(x,y,z)=1,\;z<0,\;2x>z\}$.
\end{enumerate}
\ep
\pf
Consider the diffeomorphism
\[\Phi:\,\{h\neq 0\}\cap\{z=1\}\to \{h=1\}\cap \{z\neq 0\},\quad (x,y,1)\mapsto \frac{1}{\sqrt[3]{h(x,y,1)}}(x,y,1)\]
with inverse $\Phi^{-1}(\wt{x},\wt{y},\wt{z})=\left(\frac{\wt{x}}{\wt{z}},\frac{\wt{y}}{\wt{z}},1\right).$ The restriction of $\Phi$ gives diffeomorphisms
\begin{align*}
V_+:=\{h>0\}\cap\{z=1\}&\stackrel{\approx}{\longrightarrow} \{h=1\}\cap \{z>0\}=:\mathcal{H}_+,\\
V_-:=\{h<0\}\cap\{z=1\}&\stackrel{\approx}{\longrightarrow} \{h=1\}\cap \{z<0\}=:\mathcal{H}_-.
\end{align*}
The discriminant $a^3-27b^2$ determines the number of real roots of \[f(x):=y^2-h(x,y,1)=4x^3-ax-b.\]
\begin{enumerate}
\item[a)] Case $a^3-27b^2<0$: For negative discriminant, $f(x)=4x^3-ax-b$ has exactly one real root that we denote by $x_1$. Then $V_+=\{(x,y,1)\,|\,x\in \mathbb{R},~y^2>f(x)\}$ and $V_-=\{(x,y,1)\,|\,x> x_1,~y^2<f(x)\}$ are connected. Hence, $\mathcal{H}_+=\Phi(V_+)$ and $\mathcal{H}_-=\Phi(V_-)$ are connected as well.
With $\mathcal{H}_0:=\{h=1\}\cap \{z=0\}=\{(-\frac{1}{\sqrt[3]{4}},y,0)\,|\,y\in\mathbb{R}\}$, we have $\{(x,y,z)\in\mathbb{R}^3\,|\,h(x,y,z)=1\}=\mathcal{H}_+\cup \mathcal{H}_0\cup \mathcal{H}_-$.
For $x<x_1$, we have
\[\mathcal{H}_+\ni \Phi(x,0,1)=\frac{1}{\sqrt[3]{-4x^3+ax+b}}(x,0,1)\stackrel{{x\to-\infty}}{\longrightarrow}(-\frac{1}{\sqrt[3]{4}},0,0)\in \mathcal{H}_0 \]
and for $x>x_1$,
\[\mathcal{H}_-\ni \Phi(x,0,1)=\frac{1}{\sqrt[3]{-4x^3+ax+b}}(x,0,1)\stackrel{{x\to+\infty}}{\longrightarrow}(-\frac{1}{\sqrt[3]{4}},0,0)\in \mathcal{H}_0. \]
Thus, $\{(x,y,z)\in\mathbb{R}^3\,|\,h(x,y,z)=1\}$ is connected.

\item[b)] Case $a^3-27b^2>0$: Without loss of generality, we set $a=3$ (see Lemma \ref{lemmaEquivalencePolynomial}). Then $-1<b<1$. For positive discriminant, $f(x)=4x^3-3x-b$ has three real roots that we denote by $x_1,x_2,x_3$ such that $x_2<x_3<x_1$. Hence, $f(x)>0$ for $x\in(x_2,x_3)\cup (x_1,\infty)$. As before, $V_+=\{(x,y,1)|x\in \mathbb{R},\;y^2>f(x)\}$ is connected and diffeomorphic to $\mathcal{H}_+$. $V_-$ has two connected components given by \[V_{-b}:=\{(x,y,1)|x_2<x<x_3,\;y^2<f(x)\}\text{ and } V_{-u}:=\{(x,y,1)|x> x_1,\;y^2<f(x)\}.\]
$\mathcal{H}_{-b}:=\Phi(V_{-b}) \approx V_{-b}$ and $\mathcal{H}_{-u}:=\Phi(V_{-u})\approx V_{-u}$
are connected and with $\mathcal{H}_0:=\{h=1\}\cap \{z=0\}=\{(-\frac{1}{\sqrt[3]{4}},y,0)|y\in\mathbb{R}\}$, we have
$\{(x,y,z)\in\mathbb{R}^3|h(x,y,z)=1\}=\mathcal{H}_+\cup \mathcal{H}_0\cup\mathcal{H}_{-b}\cup \mathcal{H}_{-u}$. By the same reasoning as in the proof of $a)$, we see that $\mathcal{H}_+\cup \mathcal{H}_0\cup \mathcal{H}_{-u}$ is connected.

Notice that the minimum of $f(x)$ is located at $x=\frac{1}{2}$, so $x_3<\frac{1}{2}<x_1$. Hence, we have $\mathcal{H}_{-b}=\Phi(V_{-b})\subset \{z<0\}\cap \{\frac{x}{z}<\frac{1}{2}\}$ and $\mathcal{H}_{-u}=\Phi(V_{-u}) \subset \{z<0\}\cap \{\frac{x}{z}>\frac{1}{2}\}$. From $\overline{\{z<0\}\cap\{2x> z\}}\cap\mathcal{H}_0=\emptyset$, it follows that $\overline{\mathcal{H}_{-b}}\cap (\mathcal{H}_+\cup \mathcal{H}_0 \cup \mathcal{H}_{-u})=\emptyset $. Thus, $\{(x,y,z)\in\mathbb{R}^3|h(x,y,z)=1\}$ has two connected components, namely $\mathcal{H}_+\cup \mathcal{H}_0 \cup \mathcal{H}_{-u}$ and $\mathcal{H}^{(3,b)}:=\mathcal{H}_{-b}\subset\{z<0\}\cap \{2x>z\}$.
\qed
\end{enumerate}

Now, we show that $h$ is hyperbolic for each point in the closed surface $\mathcal{H}^{(3,b)}=\{(x,y,z)\in\mathbb{R}^3| h^{(3,b)}(x,y,z)=1,\;z<0,\;2x>z\}\subset\mathbb{R}^3$ defined in Proposition \ref{propConnComponents} b).

\bp \label{propCandidate}
Let $h=y^2z-4x^3+3xz^2+bz^3$ with $-1<b<1$ and let $\mathcal{H}^{(3,b)}=\{(x,y,z)\in\mathbb{R}^3| h(x,y,z)=1,\;z<0,\,2x>z\}$. Then $-\partial^2 h(p)$ has Lorentzian signature for all $p\in \mathcal{H}^{(3,b)}$.
\ep
\pf
Let $x_2<x_3<x_1$ be the three distinct solutions of $4x^3-3x-b=0$ and let $d:=\text{det}\,\partial^2 h=-24(12xz(x+bz)+3z^3-4xy^2)$. Let 
$(x,y,1)\in V_{-b}:=\{(x,y,1)\,|\,x_2<x<x_3,~y^2<4x^3-3x-b\}$.
We show that $d(x,y,1)<0$:
\begin{description}
\item[Case $x=0$:] $-\frac{d(0,y,1)}{24}=3>0.$
\item[Case $x<0$:]
$-\frac{d(x,y,1)}{24}=12x(x+b)+3-4xy^2=12(x+\frac{b}{2})^2+3(1-b^2)-4xy^2>0.$
\item[Case $x>0$:] We use $y^2<4x^3-3x-b$:
\begin{align*}
-\frac{d(x,y,1)}{24}&=12x(x+b)+3-4xy^2>12x(x+b)+3-4x(4x^3-3x-b)\\
&=-16x^4+24x^2+16bx+3=:g(x).
\end{align*}
$g'(x)=-16(4x^3-3x-b)$, so $g(x)$ has local maxima at $x_2<-\frac{1}{2}$ and $x_1>\frac{1}{2}$ and a local minimum at $x_3$.
For $-1<b<1$, the quartic polynomial $g(x)$ has only two real roots\footnote{This can be shown using the relation between the roots of a quartic polynomial and the roots of its {\it cubic resolvent} (see e.g. \cite{I}, section 10.5). The cubic resolvent $s(z)=z^3-3z^2+3z-b^2=(z-1)^3+3-b^2$ of $g(x)$ has only one real root. This implies that $g(x)$ has two real roots.}, which must lie in $(-\infty,x_2)$ and $(x_1,\infty)$. So $g(x_3)>0$ and it follows that $g(x)>0$ for $0<x<x_3$.
\end{description}
(Note that, depending on the value of $b$, some of these cases might be empty.)

We have $\mathcal{H}^{(3,b)}=\{\frac{1}{\sqrt[3]{h(x,y,1)}}(x,y,1)\,|\,(x,y,1)\in V_{-b}\}$, where $h|_{V_{-b}}<0$ (see the proof of Proposition \ref{propConnComponents}.b). Since $d$ is homogeneous of degree three, $d|_{V_{-b}}<0$ implies $d|_\mathcal{H}>0$. Now $\partial^2 h|_{\mathcal{H}}$ can have signature $(+,+,+)$ or $(+,-,-)$. Since $\partial^2h(\partial_y,\partial_y)=2z<0$, the signature of $\partial^2 h|_{\mathcal{H}}$ is $(+,-,-)$.
\qed

\bl \label{lemmaLorentzian}
For $h^{(a,b)}=y^2z-4x^3+axz^2+bz^3$, $-\partial^2 h^{(a,b)}(x,y,z)$ is Lorentzian iff
\[(x,y,z)\in \mathcal{U}_{Lor.}:=\{\det\partial^2 h>0,\;x> 0\}\cup \{\det\partial^2 h>0,\;z< 0\}.\]
\el
\pf
For $h=y^2z-4x^3+axz^2+bz^3$, $\partial^2 h$ is positive definite iff $h_{xx}=-24x>0$, $\det\begin{pmatrix} h_{xx} & h_{xy} \\ h_{yx} & h_{yy} \end{pmatrix}=-48xz>0$ and $\det\partial^2 h=-8(12xz(ax+3bz)+a^2z^3-12xy^2)>0$, i.e.
\[\{(x,y,z)\in\mathbb{R}^3|\partial^2 h(x,y,z)>0\}=\{(x,y,z)\in\mathbb{R}^3|x<0,\;z>0,\;\det\partial^2 h>0\}=:\mathcal{U}_{pos.}.\]
Hence, $-\partial^2 h(x,y,z)$ is Lorentzian iff
\[(x,y,z)\in \{(x,y,z)\in\mathbb{R}^3\,|\,\det\partial^2 h>0\}\backslash U_{pos.}=\{\det\partial^2 h>0,~x\geq 0\}\cup \{\det\partial^2 h>0,~z\leq 0\}.\]
We show that $\{\det\partial^2 h>0,~x=0,~z\geq 0\}\cup \{\det\partial^2 h>0,~x\leq 0,~z=0\}=\emptyset:$
\begin{align*}
x=0,~z\geq 0 &\Rightarrow \det\partial^2 h=-8a^2z^3\leq 0,\\
x\leq 0,~z=0 &\Rightarrow \det\partial^2 h=96xy^2\leq 0.
\end{align*}\qed

Using the above lemma, we show that except for $\mathcal{H}^{(3,b)}$, all connected components of $\{x\in\mathbb{R}|h^{(3,b)}(x)=1,~g_x>0\}$ have nonempty boundary in $\mathbb{R}^3$.

\bp \label{propBoundary}
Let
\begin{enumerate}
\item[a)] $\mathcal{S}:=\{(x,y,z)\in\mathbb{R}^3| h(x,y,z)=1\}$ for $h=y^2z-4x^3+axz^2+bz^3$ with negative discriminant or
\item[b)] $\mathcal{S}:=\{(x,y,z)\in\mathbb{R}^3| h(x,y,z)=1\}\backslash\mathcal{H}^{(3,b)}$ for $h=y^2z-4x^3+3xz^2+bz^3$ with $-1<b<1$,
where $\mathcal{H}^{(3,b)}=\{(x,y,z)\in\mathbb{R}^3| h(x,y,z)=1,\;z<0,\;2x>z\}$.
\end{enumerate}
Then $\mathcal{S}\cap \{(x,y,z)\in\mathbb{R}^3|-\partial^2 h(x,y,z) \text{ Lorentzian}\}$ has no connected component that is closed in $\mathcal{S}$.
\ep
\pf
Let $\widetilde{\mathcal{H}}$ be a connected component of $\mathcal{S}\cap \mathcal{U}_{Lor.}$ (see lemma \ref{lemmaLorentzian}). $(-\frac{1}{\sqrt[3]{4}},0,0)\in \mathcal{S}\cap \{\det\partial^2 h=0\}$ implies $\widetilde{\mathcal{H}}\neq \mathcal{S}$. $\mathcal{U}_{Lor.}=\{\det\partial^2 h>0,~x> 0\}\cup \{\det\partial^2 h>0,~z< 0\}$ is open in $\mathbb{R}^3$, so $\mathcal{S}\cap \mathcal{U}_{Lor.}$ and hence $\widetilde{\mathcal{H}}$ are open in $\mathcal{S}$. According to proposition \ref{propConnComponents}, $\mathcal{S}$ is connected. Since $\widetilde{\mathcal{H}}\subset \mathcal{S}$ is open, nonempty and $\neq \mathcal{S}$, it cannot be closed in $\mathcal{S}$.
\qed

In summary, we have proven the following
\bc\label{corWei}
Up to linear equivalence, the only closed (in $\mathbb{R}^3$) projective special real surfaces defined by Weierstra\ss ~cubic polynomials with nonzero discriminant are given by
\[\mathcal{H}^{(3,b)}=\{(x,y,z)\in\mathbb{R}^3| y^2z-4x^3+3xz^2+bz^3=1,\;z<0,\;2x>z\}\]
with $-1<b<1$.\qed
\ec

To prove the completeness of the projective special real metric defined on the closed component $\mathcal{H}^{(3,b)}$ given in Corollary \ref{corWei}, we will use the following three lemmata.

\bl\label{lemmaR}
For each $b\in(-1,1)$ there exists exactly one $R\in\bR$, such that $h^{(3,b)}=y^2z-4x^3+3xz^2+bz^3$ is equivalent to
\[
h(x,y,z):=y^2z-x^3+xz^2+Rx^2z.
\]
\
\el
\pf It is straightforward to see that $h^{(3,b)}$ is equivalent to $h_1(x,y,z):=y^2z-x^3+xz^2+\frac{2b}{3^{\frac{3}{2}}}z^3$.
%the corresponding linear transformation is given by
%\[
%\left(
%\begin{matrix}
%\frac{1}{2^{\frac{2}{3}}} & 0 & 0\\
%0 & \frac{3^\frac{1}{4}}{2^\frac{1}{6}} & 0\\
%0 & 0 & \frac{2^\frac{1}{3}}{3^\frac{1}{2}}
%\end{matrix}
%\right).
%\]
To eliminate the $z^3$-part we make the ansatz $x=\widetilde{x}+cz$. We obtain 
\[
h_1(\widetilde{x}+cz,y,z)=y^2z-\widetilde{x}^3+(1-3c^2)\widetilde{x}z^2-3c\widetilde{x}^2z+ \left(-c^3+c+\frac{2b}{3^\frac{3}{2}}\right)z^3.
\]
We need to analyse the solvability of $\left(-c^3+c+\frac{2b}{3^\frac{3}{2}}\right)=0$. Therefore, we define $f(c):=c^3-c$ and calculate its first derivative $f'(c)=3c^2-1$. $f'(c)$ vanishes if and only if $c=-\frac{1}{3^\frac{1}{2}}$ or $c=\frac{1}{3^\frac{1}{2}}$, and $f\left(-\frac{1}{3^\frac{1}{2}}\right)=\frac{2}{3^\frac{3}{2}},\ \ \ f\left(\frac{1}{3^\frac{1}{2}}\right)=-\frac{2}{3^\frac{3}{2}}$.
Thus, for each $b\in(-1,1)$ there is exactly one $c_b\in \left(-\frac{1}{3^\frac{1}{2}},\frac{1}{3^\frac{1}{2}}\right)$, such that $f(c_b)=\frac{2b}{3^\frac{3}{2}}$. Choosing $c$ that way, we arrive at 
\[
h_1(x,y,z)=y^2z-\widetilde{x}^3+(1-3c^2)\widetilde{x}z^2-3c\widetilde{x}^2z.
\]
It follows from $|c|<\frac{1}{3^\frac{1}{2}}$ that $(1-3c^2)>0$, regardless of the choice of the parameter $b\in (-1,1)$, and, hence, the transformation $z=(1-3c^2)^{-\frac{1}{2}}\widetilde{z}$ does not switch signs. After the additional transformation $y=(1-3c^2)^{\frac{1}{4}}\widetilde{y}$, $h_1$ reads
\[
h_1(x,y,z)=\widetilde{y}^2\widetilde{z}-\widetilde{x}^3+\widetilde{x}\widetilde{z}^2+\frac{-3c}{(1-3c^2)^{\frac{1}{2}}}\widetilde{x}^2\widetilde{z}.
\]
One can easily verify that $R:\left(-\frac{1}{3^\frac{1}{2}},\frac{1}{3^\frac{1}{2}}\right)\to \bR,c\mapsto \frac{-3c}{(1-3c^2)^{\frac{1}{2}}}$ is a bijection. \qed

In the following we will omit the tildes so that $x,y,z$ denote our new coordinates.
\bl\label{lemmaWeiEquiv}
The closed surface $\mathcal{H}^{(3,b)}=\{(x,y,z)\in\mathbb{R}^3| h^{(3,b)}=1,\;z<0,\;2x>z\}$ is linearly equivalent to
\[
\mathcal{H}:=\{(x,y,z)\in\bR^3|h(x,y,z)=1,\ x>0,\ z<0\},
\]
where $h=y^2z-x^3+xz^2+Rx^2z$ as in the above lemma.
\el
\pf After transforming $h^{(3,b)}$ into $h(x,y,z)=y^2z-x^3+xz^2+Rx^2z$, the corresponding inequalities for the coordinates introduced in 
the proof of Lemma \ref{lemmaR} read
\[
z<0,\ \ x>\frac{1}{\sqrt{3}}\sqrt{\frac{1-\sqrt{3}c}{1+\sqrt{3}c}}\,z,
\]
where $c$ is the parameter determined in the proof of  Lemma \ref{lemmaR}, which satisfies $0\leq |c|<3^{-\frac{1}{2}}$. 
We have to show that $x>0$. Since $\mathcal{H}^{(3,b)}$ is connected, $\widetilde{\mathcal{H}}:=\{(x,y,z)\in\bR^3|h(x,y,z)=1,\ z<0,x>\frac{1}{\sqrt{3}}\sqrt{\frac{1-\sqrt{3}c}{1+\sqrt{3}c}}z\}$ is connected as well. Hence, it suffices to show that $\widetilde{\mathcal{H}}\cap \{x=0\}=\emptyset$. To do so, we write $\mathcal{H}$ as a graph.
\begin{align}
&h(x,y,z)=1\nonumber\\
&\Leftrightarrow y^2z-x^3+xz^2+Rx^2z=1\nonumber\\
&\Leftrightarrow z^2+\left(\frac{y^2+Rx^2}{x}\right)z-\left( x^2+\frac{1}{x}\right)=0\nonumber\\
&\Leftrightarrow z(x,y):=z=-\frac{1}{2x}\left(y^2+Rx^2+\sqrt{(y^2+Rx^2)^2+4x^4+4x}\right)\label{paraZ}.
\end{align}

The last equivalence holds since $z<0$. 

Considering the limit $x\ra 0$, we see that in the case $y\ne 0$, $\lim_{x\to 0, x>0}z(x,y)=-\infty$. 
In the case $y=0$, it is easy to verify that
\[
|z(x,0)|\geq \frac{1}{2x}\left(\sqrt{R^2x^4+4x}-|R|x^2\right)=\frac{1}{2}\sqrt{R^2x^2+\frac{4}{x}}-\frac{1}{2}|R|x\geq 0.
\]
Therefore, we have 
\[
\lim_{x\to 0, x>0}|z(x,0)|=\lim_{x\to 0, x>0}\sqrt{\frac{1}{x}}=\infty.
\]
Since $(\overline{\mathcal{H}}\setminus\{z=0\})\cap\{x=0\}=\widetilde{\mathcal{H}}\cap\{x=0\}$, the above calculation shows that $\widetilde{\mathcal{H}}\cap\{x=0\}=\emptyset$. Hence, $\widetilde{\mathcal{H}}$ can be written as the disjoint union of $\mathcal{H}$ and 
 $\{h(x,y,z)=1|z<0, 0>x>\frac{1}{\sqrt{3}}\sqrt{\frac{1-\sqrt{3}c}{1+\sqrt{3}c}}z\}$. This shows that $\mathcal{H}$ is a connected component of $\widetilde{\mathcal{H}}$. Since $\widetilde{\mathcal{H}}$ is connected, the equality $\mathcal{H}=\widetilde{\mathcal{H}}$ holds true.\qed

\bl\label{lemmaCompleteness}
Let $(M,g_1)$ be a complete Riemannian manifold and let $g_{2,p}$ be a family of complete Riemannian metrics on $\bR$, depending smoothly on $p\in M$ such that $g_{2,p}=G(p)ds^2$. Then $(M\times\bR,g_1+g_{2,p})$ is also a complete Riemannian manifold.
\el
\pf This is a special case of Theorem 2 in \cite{CHM}. \qed

\bp \label{propWeiCompl}
The surface $\mathcal{H}$ defined in Lemma \ref{lemmaWeiEquiv} endowed with the metric $g_\mathcal{H}:=-\partial^2h|_{T\mathcal{H}\times T\mathcal{H}}$ is a complete Riemannian manifold for all $R\in\bR$.
\ep
\pf Computing the Hessian of $h$, we obtain
\[
-\partial^2h=(6x-2Rz)dx^2-2zdydy-2xdzdz-4(z+Rx)dxdz-4ydydz.
\]
It was shown in the proof of Lemma 2 that $\mathcal{H}$ admits the following parametrization:
\[
F:\bR^{>0}\times\bR\to \mathcal{H}\subset\bR^3,\ (x,y)\mapsto (x,y,z(x,y)),
\]
where $z(x,y)$ is the function defined in equation $\eqref{paraZ}$.
With the abbreviation $A:=(y^2+Rx^2)^2+4x^4+4x$, $g_\mathcal{H}$ reads
\[
\begin{aligned}
g_\mathcal{H}&=x^{-3}(A^\frac{1}{2}(ydx-xdy)^2+A^{-1}((ydx-xdy)^2(12x^4y^2+6xy^2+3R^2x^4y^2+3Rx^2y^4\\
&+R^3x^6+4Rx^6+y^6)+6x^2((1+R^2x^3+4x^3)dx^2+xy^2dy^2)+2Rx^3(y^2dx^2+x^2dy^2\\
&+(ydx+xdy)^2))).
\end{aligned}
\]
Now we change the coordinates via
\[
T:\bR^{>0}\times\bR\to \bR^{>0}\times\bR,\ (s,t)\mapsto (s,st),
\]
such that, in particular, $T^*(ydx-xdy)=-s^2dt$. In these coordinates, $g_\mathcal{H}$ has the following form:
\[
\begin{aligned}
g_\mathcal{H}&=(sA)^{-1}((24s^3+6+6s^3(t^2+R)^2) ds^2\\
&+12s^4t(t^2+R)dsdt\\
&+\Big(4s^8R+s^8R^3+4Rs^5+12s^8t^2+s^8t^6+12s^5t^2+3R^2s^8t^2+ 3Rs^8t^4\\
&+(s^6R^2+2s^6t^2R+s^6t^4+4s^6+4s^3)\sqrt{s^4(t^2+R)^2+4s^4+4s}
\Big)dt^2),
\end{aligned}
\]
where $A=s^4(t^2+R)^2+4s^4+4s$ is the same function as above in the new coordinates $s$ and $t$. To show that $g_\mathcal{H}$ is complete, we start with rewriting it and make some estimates: \[
\begin{aligned}
g_\mathcal{H}&=(sA)^{-1}((24s^3+6+\frac{3}{2}s^3(t^2+R)^2) ds^2\\
&+\Big(4s^8R+s^8R^3+4Rs^5+12s^8t^2+s^8t^6+4s^5t^2+3R^2s^8t^2+ 3Rs^8t^4\\
&+(s^6R^2+2s^6t^2R+s^6t^4+4s^6+4s^3)\sqrt{s^4(t^2+R)^2+4s^4+4s}
\Big)dt^2\\
&+\underbrace{\frac{9}{2}s^3(t^2+R)^2ds^2+12s^4t(t^2+R)dsdt+8s^5t^2dt^2}_{ (\frac{3}{\sqrt{2}}s^\frac{3}{2}(t^2+R)ds+\sqrt{8}s^\frac{5}{2}tdt)^2\geq 0})\\
&\geq (sA)^{-1}((4s^3+4+s^3(t^2+R)^2) ds^2\\
&+\Big(4s^8R+s^8R^3+4Rs^5+4s^8t^2+s^8t^6+4s^5t^2+3R^2s^8t^2+ 3Rs^8t^4\\
&+(s^6R^2+2s^6t^2R+s^6t^4+4s^6+4s^3)\sqrt{s^4(t^2+R)^2+4s^4+4s}
\Big)dt^2)\\
&=\frac{1}{s^2}ds^2
+(sA)^{-1}s^2A(s^2R+s^2t^2+\sqrt{s^4(t^2+R)^2+4s^4+4s})dt^2\\
%&=\frac{1}{s^2}ds^2\\
%&+s(s^2R+s^2t^2+\sqrt{s^4(t^2+R)^2+4s^4+4s})dt^2\\
&\geq \frac{1}{s^2}ds^2
+s(s^2R+s^2t^2+\sqrt{s^4(t^2+R)^2+s^4})dt^2\\
&=\frac{1}{s^2}ds^2
+s^3(t^2+R^2+\sqrt{(t^2+R)^2+1})dt^2.
\end{aligned}
\]
Solving the ODE  
\be
\mu'=\sqrt{t^2+R+\sqrt{(t^2+R)^2+1}}, \label{muEqu} 
\ee
we obtain a  (strictly increasing) diffeomorphism $\mu : \bR \ra \bR$, $t\mapsto  \mu(t)$, such that
$d\mu^2=(t^2+R+\sqrt{(t^2+R)^2+1})dt^2$.  In fact, the right-hand side of 
\re{muEqu} is bounded from below by a positive constant.  
Now 
we can conclude the proof using Lemma \ref{lemmaCompleteness}, 
which shows that the metric $\frac{1}{s^2}ds^2
+s^3(t^2+R^2+\sqrt{(t^2+R)^2+1})dt^2=d\s^2 + e^{3\s}d\mu^2$ is complete, 
where $\s = \ln s$. 
 \qed

\section{Curvature formulas for  the generalized supergra\-vity r-map}\label{Curv_genSec}
We generalize the definition of projective special real manifolds in section \ref{sectionStatement} to hypersurfaces $\mathcal{H}\subset\mathbb{R}^{n+1}$ defined as the level set of an arbitrary homogeneous function $h$ and then define the generalized supergravity r-map, which assigns to each such $n$-dimensional real manifold a K\"ahler manifold of dimension $2m:=2n+2$ (see \cite{CHM}). Then we calculate the Riemann, Ricci and scalar curvature of manifolds in the image of the generalized supergravity r-map and express them in terms of $h$ and its derivatives.

Let $U\subset \mathbb{R}^{n+1}$ be an $\mathbb{R}^{>0}$-invariant domain and let $h:U\to \mathbb{R}^{>0}$ be a smooth function which is homogeneous of degree\footnote{Note that one can extend this definition to $D<0$. For $D<0$, we then get a different signature for $g_M$ compared to the conventions in \cite{CHM}.} $D\in \mathbb{R}^{>0}\backslash \{1\}$. Then
$\mathcal{H}:=\{x\in U|h(x)=1\}\subset U$ is a smooth hypersurface and we will assume that $-\partial^2h$ induces a Riemannian metric $g_{\mathcal{H}}=-\partial^2h|_{T\mathcal{H}\times T\mathcal{H}}$ on $\mathcal{H}$.

We consider $M:=\mathbb{R}^{n+1}+iU\subset \mathbb{C}^{n+1}$, endowed with the standard complex structure and the standard holomorphic coordinates $z^1,\ldots, z^{n+1}$ induced from $\mathbb{C}^{n+1}$. We define $x^\mu:=\text{Im}\,z^\mu$ and $y^\mu:=\text{Re}\,z^\mu$, i.e. $z=(z^\mu)_{\mu=1,\ldots,n+1} =y+ix$ with $x=(x^\mu)\in U$ and $y=(y^\mu)\in \mathbb{R}^{n+1}$. We define a positive definite K\"ahler metric\footnote{ Note that $g_M$ differs from the metric in \cite{CHM} by a conventional factor: $g_M^{[CHM]}=\frac{3}{2D}g_M$.}
$g_M$ on $M$ with K\"ahler potential $K(z,\bar{z})=-\log h(x)$:
\begin{equation}\label{rMapMetric}
g_M=2g_{\mu\bar{\nu}}dz^\mu d\bar{z}^\nu,~ g_{\mu\bar{\nu}}:=\frac{\partial^2K}{\partial z^\mu \partial\bar{z}^\nu}=-\frac{h_{\mu\nu}(x)}{4h(x)}+\frac{h_\mu(x)h_\nu(x)}{4h^2(x)}.
\end{equation}
Here, we use the notation $h_\mu(x):=\frac{\partial{h}(x)}{\partial x^\mu}$, $h_{\mu\nu}(x):=\frac{\partial^2{h}(x)}{\partial x^\mu \partial x^\nu}$, \ldots.

\bd
We call the correspondence $(\mathcal{H},\,g_{\mathcal{H}})\mapsto (M,\,g_M)$ the {\cmssl generalized supergravity r-map}. The restriction to polynomial functions $h$ of degree $D=3$ is called the {\cmssl supergravity r-map}. Manifolds in the image of the supergravity r-map are called {\cmssl projective very special K\"ahler}.
\ed

Note that the K\"ahler manifolds in the image of the generalized supergravity r-map in general only fall into the class of projective special K\"ahler manifolds (see e.g. \cite{CHM}) if $h$ is a polynomial of degree $D=3$.

Using the fact that $h$ is a homogeneous function of degree $D$, i.e.
$\sum_\mu h_\mu(x) x^\mu=D\cdot h(x)$, $\sum_\nu h_{\mu\nu}(x) x^\nu=(D-1)\cdot h_\mu(x)$, \ldots\,,
one can check that the coefficients of the inverse metric $g_M^{-1}$ are given by\footnote{This has been found by specializing the formula for the inverse of projective special K\"ahler metrics in \cite{C-G} to projective very special K\"ahler metrics ($D=3$) and then generalizing this to metrics of the form \eqref{rMapMetric} defined by arbitrary homogeneous functions $h$.}
\begin{equation}\label{InverseMetric} g^{\bar{\nu}\lambda}=-4h(x) h^{\nu\lambda}(x)+\frac{4}{D-1}x^\nu x^\lambda,
\end{equation}
where $h^{\mu\nu}$ denote the coefficients of the matrix $(h_{\mu\nu})^{-1}$.

\bt\label{curvatureTheorem}
Let $M$ be a $2m$-dimensional manifold in the image of the generalized supergravity r-map described by a homogeneous function $h$ of degree $D\in\mathbb{R}^{>0}\backslash\{1\}$. Then the Riemann, Ricci and scalar curvature\footnote{We define the scalar curvature $scal$ for K\"ahler manifolds to be one half of the trace of $Ric$, i.e. $scal:=g^{\mu\bar{\nu}}Ric_{\mu\bar{\nu}}$. Compared to the standard definition $scal_\mathbb{R}:=\tr\!\!Ric$ of the scalar curvature in differential geometric literature, we thus have $scal_\mathbb{R}=2scal$.} in the holomorphic coordinates $(z^\mu)=(y^\mu+ix^\mu)$ defined above are given by
\begin{align} R^\rho_{~\sigma\mu\bar{\nu}}&=-\frac{1}{4h^2}\Bigg[-h^2h^\rho_{~\sigma\mu\nu}+\frac{1}{D-1}x^\rho(h h_{\sigma\mu\nu}-h_{\mu\sigma}h_\nu)+h_\mu h_\nu\delta_\sigma^\rho+h_\sigma h_\nu\delta_\mu^\rho \nonumber \\
&\quad\quad -h\left( h_{\sigma\nu}\delta^\rho_\mu+h_{\mu\nu}\delta_\sigma^\rho-\frac{1}{D-1}h_{\mu\sigma}\delta_\nu^\rho\right)-h^2h^\rho_{~\nu\beta}h^\beta_{~\sigma\mu} \Bigg],
\end{align}
\begin{equation}\label{RicciCurvature} Ric_{\mu\bar{\nu}}=-mg_{\mu\bar{\nu}}+\frac{1}{4}h_{\mu}^{~\alpha\beta}h_{\alpha\beta\nu}+\frac{1}{4}h^\rho_{~\rho\mu\nu},
\end{equation}
\begin{equation} scal=-m^2+\frac{D-2}{D-1}m+hh_{\alpha\beta\gamma}h^{\alpha\beta\gamma}+hh^{\mu\nu}_{~~\mu\nu}, \end{equation}
where the argument of $h$ and its derivatives is always $x$ and where we use the Lorentzian metric $-\partial^2 h$ to raise and lower indices. Sums over repeated indices are implied via the Einstein summation convention.
\et

\pf
For K\"ahler manifolds, the only non-vanishing Christoffel symbols are (see e.g. \cite{M}, section 12.2, or \cite{KN})
\begin{equation*} dz^\rho(\nabla_{\partial_{z^\sigma}}\partial_{z^\mu})=:\Gamma^\rho_{\sigma\mu}=g^{\rho\bar{\kappa}}\partial_{z^\sigma} g_{\mu\bar{\kappa}}
\end{equation*}
and its complex conjugate.
For the Riemann tensor $R(X,Y):=\nabla_X\nabla_Y Z-\nabla_Y\nabla_XZ-\nabla_{[X,Y]}Z$ in holomorphic coordinates, we have (see e.g. \cite{M}, section 12.2, or \cite{KN})
\begin{equation*}  dz^\rho\left( R(\partial_{z^\mu},\,\partial_{z^{\bar{\nu}}}) \partial_{z^\sigma}\right)=: R^\rho_{~\sigma\mu\bar{\nu}}=-\partial_{\bar{z}^\nu}\Gamma^\rho_{\sigma\mu}.
\end{equation*}
The other non-vanishing components
$R^{\bar{\rho}}_{~\bar{\sigma}\mu\bar{\nu}}$, $R^\rho_{~\sigma\bar{\mu}\nu}$ and $R^{\bar{\rho}}_{~\bar{\sigma}\bar{\mu}\nu}$ of the curvature tensor can be obtained from this via symmetry and complex conjugation.

Since the metric only depends on the imaginary part of $z=y+ix\in M$,
we have
\begin{align*} \partial_{z^\sigma}g_{\mu\bar{\kappa}}&=-\frac{i}{2}\partial_{x^\sigma}\left( -\frac{h_{\mu\kappa}}{4h}+\frac{h_\mu h_\kappa}{4h^2}\right) \\
&=-\frac{i}{2}\frac{1}{4h^3}\left( -h^2 h_{\mu\kappa\sigma}+hh_\sigma h_{\mu\kappa} +hh_\kappa h_{\mu\sigma}+h h_\mu h_{\kappa\sigma}-2h_\mu h_\kappa h_\sigma\right),\end{align*}
where the argument of $h$ and its derivatives is always $x$. This gives
\begin{align}
\Gamma^\rho_{\sigma\mu}&=g^{\rho\bar{\kappa}}\partial_{z^\sigma}g_{\mu\bar{\kappa}}=-\frac{i}{2}\left(-4hh^{\rho\kappa}+\frac{4}{D-1}x^\rho x^\kappa\right)\nonumber\\
&\quad\quad\quad\quad\quad\quad\quad\quad\cdot \frac{1}{4h^3}\left(-h^2 h_{\mu\kappa\sigma}+hh_\sigma h_{\mu\kappa} +hh_\kappa h_{\mu\sigma}+h h_\mu h_{\kappa\sigma}-2h_\mu h_\kappa h_\sigma\right) \nonumber \\
&=-\frac{i}{2h^3}\Bigg( h^3h^{\rho\kappa}h_{\kappa\mu\sigma}-h^2h_\sigma\delta^\rho_\mu-\frac{1}{D-1}h^2x^\rho h_{\mu\sigma}-h^2h_\mu\delta_\sigma^\rho+\frac{1}{D-1}2hx^\rho h_\mu h_\sigma\nonumber \\
& \quad\quad-\frac{D-2}{D-1}h^2x^\rho h_{\mu\sigma} +h x^\rho h_\sigma h_\mu +\frac{D}{D-1}h^2x^\rho h_{\mu\sigma}+h x^\rho h_\mu h_\sigma-\frac{D}{D-1}2hx^\rho h_\mu h_\sigma \Bigg)\nonumber \\
&=-\frac{i}{2h}\left(hh^{\rho\kappa}h_{\kappa\mu\sigma}-h_\sigma \delta^\rho_\mu-h_\mu\delta^\rho_\sigma +\frac{1}{D-1}x^\rho h_{\mu\sigma}\right)\nonumber,
\end{align}
where for the third equality, we used $h^{\rho\kappa}h_\kappa=\frac{1}{D-1}x^\rho$. The curvature tensor is then found to be
\begin{align}
R^\rho_{~\sigma\mu\bar{\nu}}&=-\frac{i}{2}\partial_{x^\nu}\Gamma^\rho_{\sigma\mu}=-\frac{1}{4h}\Bigg[-\frac{h_\nu}{h}\left(hh^{\rho\kappa}h_{\kappa\mu\sigma}-h_\sigma \delta_\mu^\rho-h_\mu\delta_\sigma^\rho+\frac{1}{D-1}x^\rho h_{\mu\sigma}\right) \nonumber \\
&\quad\quad +h_\nu h^{\rho\kappa}h_{\kappa\mu\sigma}-hh^{\rho \alpha}h_{\nu \alpha \beta}h^{\beta\kappa}h_{\kappa\mu\sigma}+hh^{\rho\kappa}h_{\kappa\mu\sigma\nu} \nonumber \\
&\quad\quad -h_{\sigma\nu}\delta_\mu^\rho-h_{\mu\nu}\delta_\sigma^\rho+\frac{1}{D-1}h_{\mu\sigma}\delta_\nu^\rho+\frac{1}{D-1}x^\rho h_{\mu\sigma\nu}\Bigg]\nonumber \\
&=-\frac{1}{4h^2}\Bigg[ h^2h^{\rho\kappa}h_{\kappa\mu\sigma\nu}+\frac{1}{D-1}x^\rho(h h_{\mu\sigma\nu}-h_{\mu\sigma}h_\nu)+h_\mu h_\nu\delta_\sigma^\rho+h_\sigma h_\nu\delta_\mu^\rho \nonumber \\
&\quad\quad -h\left( h_{\sigma\nu}\delta^\rho_\mu+h_{\mu\nu}\delta_\sigma^\rho-\frac{1}{D-1}h_{\mu\sigma}\delta_\nu^\rho\right)-h^2h^{\rho \alpha}h_{\nu\alpha\beta}h^{\beta\kappa}h_{\kappa\mu\sigma} \Bigg],\label{curvatureTensor}
\end{align}
where for the second equality, we used the formula $\frac{d}{d t} A^{-1}=-A^{-1}\frac{dA}{d t} A^{-1}$ for an invertible matrix $A$ that smoothly depends on a parameter $t$. Now, we contract the indices $\rho$ and $\mu$ in equation \eqref{curvatureTensor} to obtain the Ricci tensor:
\begin{align}
Ric_{\mu\bar{\nu}}&=R^\rho_{~\mu\rho\bar{\nu}}=-\frac{1}{4h^2}\Bigg[ h^2h^{\rho\kappa}h_{\kappa\rho\mu\nu}+\frac{1}{D-1}x^\rho(h h_{\rho\mu\nu}-h_{\rho\mu}h_\nu)+h_\rho h_\nu\delta_\mu^\rho+h_\mu h_\nu\delta_\rho^\rho \nonumber \\
&\quad\quad -h\left( h_{\mu\nu}\delta^\rho_\rho+h_{\rho\nu}\delta_\mu^\rho-\frac{1}{D-1}h_{\rho\mu}\delta_\nu^\rho\right)-h^2h^{\rho \alpha}h_{\nu\alpha\beta}h^{\beta\kappa}h_{\kappa\rho\mu} \Bigg] \nonumber \\
&=-\frac{1}{4h^2}\Bigg[h^2h^{\rho\kappa}h_{\kappa\rho\mu\nu}+\frac{1}{D-1}\left((D-2)h h_{\mu\nu}-(D-1)h_\mu h_\nu\right) + h_\mu h_\nu+mh_\mu h_\nu \nonumber \\
&\quad\quad -h\left(m h_{\mu\nu}+ h_{\mu\nu}-\frac{1}{D-1}h_{\mu\nu}\right)-h^2h^{\rho\alpha}h_{\nu\alpha\beta}h^{\beta\kappa}h_{\kappa\rho\mu}\Bigg] \nonumber \\
&=-mg_{\mu\bar{\nu}}+\frac{1}{4}h^{\rho\alpha}h_{\nu\alpha\beta}h^{\beta\kappa}h_{\kappa\rho\mu}-\frac{1}{4}h^{\rho\kappa}h_{\kappa\rho\mu\nu}.\nonumber
\end{align}
The scalar curvature for manifolds in the image of the generalized local r-map reads
\begin{align}
scal&=g^{\mu\bar{\nu}}Ric_{\mu\bar{\nu}}=(-4hh^{\mu\nu}+\frac{4}{D-1}x^\mu x^\nu)\left(-mg_{\mu\bar{\nu}}+\frac{1}{4}h^{\rho\alpha}h_{\nu\alpha\beta}h^{\beta\kappa}h_{\kappa\rho\mu}-\frac{1}{4}h^{\rho\kappa}h_{\kappa\rho\mu\nu}\right)\nonumber\\
&=-m^2-hh^{\rho\alpha}h_{\nu\alpha\beta}h^{\beta\kappa}h^{\mu\nu}h_{\kappa\rho\mu}+\frac{(D-2)^2}{D-1}\delta_\beta^\rho \delta_\rho^\beta+hh^{\rho\kappa}h_{\kappa\rho\mu\nu}h^{\mu\nu}-\frac{(D-3)(D-2)}{D-1}\delta_\rho^\rho \nonumber\\
&=-m^2+\frac{D-2}{D-1}m-hh^{\rho\alpha}h_{\nu\alpha\beta}h^{\beta\kappa}h^{\mu\nu}h_{\kappa\rho\mu}+hh^{\rho\kappa}h_{\kappa\rho\mu\nu}h^{\mu\nu}.\nonumber
\end{align}
\qed

\br
\rm{Note that in the derivation of these formulas, we did not use the fact that $D$ is positive, i.e. the theorem also holds for metrics of the form \eqref{rMapMetric} with $h$ being a homogeneous function of degree $D<0$.}
\er

The curvature for projective very special K\"ahler manifolds, i.e. for manifolds in the image of the generalized supergravity r-map that are defined by a homogeneous cubic polynomial $h$ can be easily obtained from Theorem \ref{curvatureTheorem} by setting $D=3$ and dropping terms with quadruple derivatives of $h$.

For an arbitrary homogeneous function $h$, we also give the following alternative expression for $Ric$ and $scal$:
\bc\label{curvatureCor}
Let $M$ be a $2m$-dimensional manifold in the image of the generalized supergravity r-map described by a homogeneous function $h$ of degree $D\in\mathbb{R}^{>0}\backslash\{1\}$. Then the Ricci tensor and the scalar curvature are given by
\begin{equation}\label{Ricci2} Ric_{\mu\bar{\nu}}=-\partial_{z^\mu}\partial_{\bar{z}^\nu}\log \frac{d}{h^m},
\end{equation}
\begin{equation}\label{scalCurvature} scal=-m^2+\frac{D-2}{D-1}m-\frac{h}{d}d^\mu_{~\mu}+\frac{h}{d^2}d^\mu d_\mu, \end{equation}
where $d(x):=\det (\partial^2 h(x))$.
\ec
\pf
Using $\frac{d}{dt} (\det A)=\det A \cdot \tr\!\!\left(A^{-1}\frac{dA}{dt}\right)$, we get
\begin{align}-\partial_{z^\mu}\partial_{\bar{z}^\nu}\log \frac{d}{h^m}&=m\cdot \partial_{z^\mu}\partial_{\bar{z}^\nu}\log h-\frac{1}{4}\frac{\partial^2}{\partial x^\mu\partial x^\nu}\log d=-mg_{\mu\bar{\nu}}-\frac{1}{4}\frac{\partial}{\partial x^\mu}\frac{d_\nu}{d}\nonumber\\
&=-mg_{\mu\bar{\nu}}-\frac{1}{4}\frac{\partial}{\partial x^\mu}(h^{\alpha\beta} h_{\nu\alpha\beta})\nonumber\\
&=-mg_{\mu\bar{\nu}}+\frac{1}{4}h^{\alpha\rho}h_{\mu\rho\sigma}h^{\sigma\beta}h_{\nu\alpha\beta}-\frac{1}{4}h^{\alpha\beta}h_{\mu\nu\alpha\beta}\nonumber\\
&\stackrel{\eqref{RicciCurvature}}{=}Ric_{\mu\bar{\nu}}.\nonumber
\end{align}
Since $h_{\mu\nu}$ are homogeneous of degree $D-2$, $d=\det (h_{\mu\nu})$ is homogeneous of degree $m(D-2)$. Using this property and \eqref{InverseMetric}, we find
\begin{align}
scal&=g^{\mu\bar{\nu}}Ric_{\mu\bar{\nu}}=-m^2+\frac{h}{d}h^{\mu\nu}d_{\mu\nu}-\frac{h}{d^2}h^{\mu\nu}d_\mu d_\nu\nonumber\\
&\quad\quad\quad\quad~\quad\quad -\frac{1}{D-1}m(D-2)(m(D-2)-1)+\frac{1}{D-1}(m(D-2))^2 \nonumber \\
&=-m^2+\frac{D-2}{D-1}m+\frac{h}{d}h^{\mu\nu}d_{\mu\nu}-\frac{h}{d^2}h^{\mu\nu}d_\mu d_\nu.\nonumber
\end{align}
\qed

\section{Classification of complete projective very special K\"ahler manifolds of complex dimension 3}\label{lastSec}
In \cite{CHM}, it was shown that the supergravity r-map maps \emph{complete} $n$-dimensional projective special real manifolds to \emph{complete} projective special K\"ahler manifolds of complex dimension $m:=n+1$. Since there is a totally geodesic embedding of a projective special real manifold into the corresponding very special K\"ahler manifold, the image of an incomplete manifold under the r-map is incomplete. From the classification of all complete projective special real surfaces in Theorem \ref{mainThm}, we thus immediately get the following corollary:
\bc
The supergravity r-map assigns to each projective special real surface given in Theorem \ref{mainThm} a complete projective special K\"ahler 3-manifold and up to isometry, any complete projective special K\"ahler 3-manifold in the image of the r-map is obtained from one of the surfaces in Theorem \ref{mainThm}.\qed
\ec

To classify all complete projective very special K\"ahler manifolds up to isometry,  we want to show that the supergravity r-map maps the list of complete surfaces in Theorem \ref{mainThm} to a list of pairwise non-isometric manifolds.

Using the formula for the scalar curvature of manifolds in the image of the supergravity r-map given in Corollary \ref{curvatureCor}, we obtain the following result:
\bp
The five complete projective special K\"ahler manifolds in the image of the supergravity r-map obtained from the examples a)-e) in Theorem \ref{mainThm} are pairwise non-isometric.
\ep
\pf
Applied to the case of projective very special K\"ahler 3-manifolds, the formula for the scalar curvature given in Corollary \ref{curvatureCor} reads
\begin{equation}\label{scalPVSK}
scal=-\frac{15}{2}+\frac{h}{d}h^{\mu\nu}d_{\mu\nu}-\frac{h}{d^2}h^{\mu\nu}d_\mu d_\nu.
\end{equation}
We use it to determine the image of $scal: M\to \mathbb{R}$ for the five projective special K\"ahler manifolds in the image of the supergravity r-map obtained from the examples a)-e) in Theorem \ref{mainThm}. $scal(M)$ is pairwise different for the examples a)-e) and hence they are non-isometric.
\begin{enumerate}
\item[a)]
Example a), the so-called STU model, is defined by
\begin{equation*} h: U\to \mathbb{R},~ (x,y,z)\mapsto xyz ~~\text{ with }~~ U=\{(x,y,z)\in\mathbb{R}^3|x>0,\;y>0,\;z>0\}.\end{equation*}
The corresponding manifold in the image of the r-map is the symmetric space $(SU(1,1)/U(1))^3$. Consequently, its scalar curvature is constant:
\begin{equation*}d:=\text{det}\, h_{\mu\nu}=2h \quad \stackrel{\eqref{Ricci2}}{\Rightarrow}\quad Ric_{\mu\bar{\nu}}=-2g_{\mu\bar{\nu}}\quad \Rightarrow\quad scal=-6,\end{equation*}
i.e. the image of the scalar curvature is $scal(M)=\{-6\}$.
\item[b)]
This example is described by
\begin{equation*} h: U\to \mathbb{R},~ (x,y,z)\mapsto x(xy-z^2) ~~\text{ with }~~ U=\{(x,y,z)\in\mathbb{R}^3|x>0,\;y>\frac{z^2}{x}\}.\end{equation*}
While the corresponding projective special real manifold is the symmetric space $SO(2,1)/SO(2)$, the corresponding projective special K\"ahler manifold obtained from the r-map is homogeneous but non-symmetric \cite{DV1}.

We have
$(h^{\mu\nu})=\frac{1}{2x^3}\begin{pmatrix} 0 & x^2 & 0 \\ x^2 & -z^2-xy & -xz \\ 0 & -xz & -x^2 \end{pmatrix}$ and with $d:=\det(h_{\mu\nu})=8x^3$, one finds $h^{\mu\nu}d_{\mu\nu}=0$, $h^{\mu\nu}d_\mu d_\nu=0$. \eqref{scalPVSK} then gives $scal(M)=\{-7.5\}$.
\item[c)]
The so-called quantum STU model is defined by
\begin{equation*} h: U\to \mathbb{R},~ (x,y,z)\mapsto xyz+x^3 ~~\text{ with }~~ U=\{(x,y,z)\in\mathbb{R}^3|x<0,~z<0,~y>-\frac{x^2}{z}\}.\end{equation*}
With $d=2h-8x^3$ and $h^{\mu\nu}=\frac{1}{d}\begin{pmatrix} -x^2 & xy & xz \\ xy & -y^2 & yz-6x^2 \\ xz & yz-6x^2 & -z^2 \end{pmatrix}$,
we calculate
\begin{equation*}
h^{\mu\nu}d_{\mu\nu}=12\cdot\frac{h}{d},\quad h^{\mu\nu}d_\mu d_\nu=\frac{12}{d}\cdot\left[ (xyz-5x^3)^2-52x^6\right].
\end{equation*}
The scalar curvature can be written as
\begin{equation*} scal=-\frac{15}{2}+3h\cdot\left(\frac{1}{d}+48x^3\cdot\frac{h}{d^3}\right).\end{equation*}
Using $d=2h-8x^3$, we can show that $\text{scal}>-\frac{15}{2}$:
\begin{equation*} \text{scal}>-\frac{15}{2}~~\stackrel{(h,d>0)}{\Leftrightarrow}~~  d^2+48x^3h>0 ~~ \Leftrightarrow ~~ (2h+4x^3)^2+48x^6>0.\end{equation*}
We also check that $\text{scal}<-6$:
\begin{align} \text{scal}<-6~~&\Leftrightarrow~~ 3h\left(\frac{1}{d}+48x^3\frac{h}{d^3}\right)<\frac{3}{2}& ~~ &\stackrel{(d>0)}{\Leftrightarrow} ~~ 2h(1+48x^3\frac{h}{d^2})<d=2h-8x^3 \nonumber \\
&\Leftrightarrow ~~ 2\cdot 48 x^3h^2<-8x^3d^2& ~~ &\stackrel{(x<0)}{\Leftrightarrow}~~ 12h^2+
d^2>0.\quad\nonumber \end{align}
To show that the bounds $-\frac{15}{2}<\text{scal}<-6$ are optimal, we determine the behaviour of $\text{scal}$ at the boundary $\partial U=\{y=-\frac{x^2}{z},~x<0,~z<0\}\,\cup\, \{x=0,~y\geq0,~z\leq0\}:$
For $\partial U\cap \{x< 0\}$, we have $h|_{\partial U\cap \{x< 0\}}=0$, $d|_{\partial U\cap \{x< 0\}}=-8x^3\neq 0$ and hence \begin{equation*}\text{scal}\stackrel[x=x_0<0]{h\to 0}{\longrightarrow}-\frac{15}{2}.\]
For $\{x= 0,~y>0,~z<0\}\subset\partial U\cap \{x= 0\}$, we have \begin{equation*}\text{scal}\stackrel[y=y_0>0,~z=z_0< 0]{ x\rightarrow 0}{\longrightarrow} \frac{-15}{2}+\frac{3}{2}=-6.\end{equation*}

Since $U$ is connected and $\text{scal}$ is continuous, we have proven that
\begin{equation*} \text{scal}(M)=\Big(-\frac{15}{2},-6\Big)=(-7.5,-6). \end{equation*}

\item[d)]
This example is described by
\begin{equation*} h: U\to \mathbb{R},~ (x,y,z)\mapsto z(x^2+y^2-z^2) ~~\text{with}~~ U=\{(x,y,z)\in\mathbb{R}^3\,|\,z<0,~x^2+y^2<z^2\}.\end{equation*}

Using
$d=-8(h+4z^3)$ and $h^{\mu\nu}=\frac{-4}{d}\begin{pmatrix} y^2+3z^2 & -xy & xz \\ -xy & x^2+3z^2 & yz \\ xz & yz & -z^2 \end{pmatrix}$,
we calculate
\begin{equation*}
h^{\mu\nu}d_{\mu\nu}=192\cdot\frac{h}{d},\quad h^{\mu\nu}d_\mu d_\nu=\-12d+4\cdot 96z^3+4\cdot96^2z^6\frac{1}{d}.
\end{equation*}

The scalar curvature can be written as
\begin{equation*} scal=-\frac{15}{2}-12h\cdot\left(\frac{1}{d}-16\cdot48z^3\cdot \frac{h}{d^3}\right).\end{equation*}
We show that $\text{scal}<-\frac{15}{2}$:
\begin{equation} \text{scal}<-\frac{15}{2}~~\stackrel{(h,d>0)}{\Leftrightarrow}~~ d^2-16\cdot48z^3h>0 ~~ \Leftrightarrow ~~ (h-2z^3)^2+12z^6>0.\quad\nonumber \end{equation}
This bound is assumed at the boundary of $U$: \[\stackrel[p\to p_0]{ }{\lim} scal_p=-\frac{15}{2}~~ \forall \, p_0\in\partial U\backslash\{0\}.\end{equation*}
At $(x,y,z)=(0,0,-1)\in U$, we have $\text{scal}(0,0,-1)=-8-\frac{2}{3}$. To show that $-8-\frac{2}{3}$ is a lower bound for the scalar curvature, we show $-\frac{12h}{d}\geq -\frac{1}{2}$ and $12\cdot 16\cdot 48 \frac{z^3h}{d^3}\geq -\frac{2}{3}$:
\begin{equation} -\frac{12h}{d}\geq -\frac{1}{2}~~\stackrel{(d>0)}{\Leftrightarrow}~~ 0\geq -d+24h ~~ \Leftrightarrow ~~ 0\geq h+z^3 ~~ \stackrel{(z<0)}{\Leftrightarrow} ~~ 0\leq x^2+y^2, \nonumber \end{equation}
\begin{align}
12\cdot 16\cdot 48 \frac{z^3h}{d^3}\geq -\frac{2}{3}~~&\Leftrightarrow~~  27z^3h^2\geq (h+4z^3)^3 \nonumber \\ ~~&\Leftrightarrow~~ 27z^5(x^2+y^2-z^2)^2\geq z^3(x^2+y^2+3z^2)^3 \nonumber \\
~~&\Leftrightarrow~~ 0\leq (x^2+y^2)^3-18(x^2+y^2)^2z^2+81(x^2+y^2)z^4 \nonumber \\
~~&\Leftrightarrow~~ 0\leq (x^2+y^2)(x^2+y^2-9z^2)^2.\nonumber \end{align}

We have thus proven that
\begin{equation*} \text{scal}(M)=\Big[-8-\frac{2}{3},-\frac{15}{2}\Big)=[-8.\bar{6},-7.5). \end{equation*}

\item[e)]
This example is described by
\begin{equation*} h: U\to \mathbb{R},~ (x,y,z)\mapsto x(y^2-z^2)+y^3 ~~\text{with}~~ U=\{(x,y,z)\in\mathbb{R}^3\,|\,y<0,~x>0,~h>0\}.\end{equation*}
One has $d=8(x(y^2-z^2)-3yz^2)$ and $h^{\mu\nu}=\frac{4}{d}\begin{pmatrix} -x(x+3y) & xy & (x+3y)z \\ xy & -z^2 & -yz \\ (x+3y)z & -yz & -y^2 \end{pmatrix}$.

Then $h^{\mu\nu}d_{\mu\nu}=192 \frac{h}{d}$
and
\begin{equation*}
h^{\mu\nu}d_\mu d_\nu = 4\cdot \frac{64}{d}\left[3x^2y^4-6x^2y^2z^2+3x^2z^4-3xy^5-24xy^3z^2+27xyz^4-72y^4z^2-9z^6\right].
\end{equation*}
The scalar curvature can be witten as
\begin{equation*} \begin{array}{l} scal=-\frac{15}{2}+192\cdot\frac{h^2}{d^2} -12\cdot 64\cdot \frac{h}{d^3}\Big[x(x-y)y^4 \\
\quad\quad\quad\quad -2y^2z^2\left((x+2y)^2+8y^2\right)+x(x+9y)z^4-3z^6\Big].\end{array} \end{equation*}
Since $scal$ only contains even powers of $z$, we have $\frac{\partial scal}{\partial z}\Big|_{z=0}=0$. We restrict ourselves to the hypersurface $M\cap \{z=0\}\subset M$ and determine $scal(M\cap \{z=0\})$.
\begin{equation*}scal(x,y,0)=-6+\frac{12xy+9y^2}{2 x^2}\end{equation*} has critical points only for $y=-\frac{2}{3}x$, where it assumes the value $-8$:
\begin{equation*}\partial_y scal(x,y,0)=\frac{6x+9y}{x^2},\quad  \partial_x scal(x,y,0)=-y\frac{6x+9y}{x^3};\quad  scal(x,-\frac{2}{3}x,0)=-8.\end{equation*}
Since $scal$ is homogeneous of degree zero, it suffices to consider the image of \[ (-1,0)\to\mathbb{R},\quad y\mapsto scal(1,y,0).\] At the boundaries $y=-1$ and $y=0$ of $U\cap\{x=1,z=0\}$, we have
\begin{equation*}\stackrel[y\to -1]{ }{\lim}scal(1,y,0)=-7.5 \quad \text{and} \quad \stackrel[y\to 0]{ }{\lim}scal(1,y,0)=-6.\end{equation*}
This shows that $scal(M\cap \{z=0\})=[-8,-6)$. In particular, we have
\begin{equation*} [-8,-6)\subset \text{scal}(M). \end{equation*}
\end{enumerate}
\qed
 
\br
\rm{Note that the results obtained in the proof of the above proposition show that the complete projective special K\"ahler manifolds obtained from examples $c)$, $d)$ and $e)$ in Theorem \ref{mainThm} via the supergravity r-map have non-constant scalar curvature and, hence, are not locally homogeneous.

Using the formula for the scalar curvature in Corollary \ref{curvatureCor}, one can similarly show that all manifolds in the one-parameter family of complete projective special K\"ahler manifolds obtained from Weierstra{\ss} cubic polynomials (see example $f)$ in Theorem \ref{mainThm}) are not locally homogeneous. This one-parameter family is particularly interesting, since using the supergravity c-map, which maps complete projective special K\"ahler manifolds to complete quaternionic K\"ahler manifolds \cite{CHM}, it gives an explicit expression for a one-parameter family of complete quaternionic K\"ahler metrics.}
\er

\br
\rm{For the one-parameter family of complete projective very special K\"ahler manifolds defined by Weierstra{\ss} polynomials (see Example f) in Theorem \ref{mainThm}), we obtain the following results using numerical methods:
\begin{equation*} scal(M)=\begin{cases} [s_{min}(b),\,s_{max}(b)] &\quad\quad\mbox{for } -1<b<0 \\
[s_{min}(b),\,-7.5) &\quad\quad\mbox{for } 0\leq b < 1, \end{cases}
\end{equation*}
where $s_{max}:\,(-1,0)\stackrel{\sim}{\to}(-7.5,-6)$ and $s_{min}:(-1,1)\stackrel{\sim}{\to}(-8.\bar{6},-8)$ are strictly decrea\-sing.
This shows that all manifolds in the image of the supergravity r-map obtained from the examples in Theorem \ref{mainThm} are non-isometric and hence it finishes the classification of all complete projective very special K\"ahler 3-manifolds.}
\er

\br
\rm{There exist precisely two complete projective special real curves, up to linear equivalence \cite{CHM}:
$\mathcal{H}_{hom.}:=\{(x,y)\in\mathbb{R}^2|x^2y=1,~x>0\}$ and $\mathcal{H}_{inh.}:=\{(x,y)\in\mathbb{R}^2|x(x^2-y^2)=1,~x>0\}$,
where $\mathcal{H}_{hom.}$ admits a transitive group of linear transformations, while $\mathcal{H}_{inh.}$ is inhomogeneous.

One can show using the curvature formulas in Theorem \ref{curvatureTheorem} that the projective special K\"ahler manifold $M_{hom.}$ obtained from $\mathcal{H}_{hom.}$ via the supergravity r-map is a product of two complex hyperbolic lines with different curvature, which is well-known from the physics literature (see e.g. \cite{DV2}).
On the other hand, the projective special K\"ahler manifold $M_{inh.}$ corresponding to $\mathcal{H}_{inh.}$ has non-constant scalar curvature and hence, it is not locally homogeneous.}
\er

\end{document}